\documentclass[12pt]{amsart}
\usepackage{amsfonts}
\usepackage{mathrsfs}
\usepackage{amscd,amsmath,amssymb,amsfonts}
\usepackage[all]{xy}
\theoremstyle{plain}
\newtheorem{thm}{Theorem}
\newtheorem{lem}[thm]{Lemma}
\newtheorem{cor}[thm]{Corollary}
\newtheorem{prop}[thm]{Proposition}

\theoremstyle{definition}
\newtheorem{defn}[thm]{Definition}

\newtheorem{question}[thm]{Question}
\newtheorem{rmk}[thm]{Remark}

\newtheorem{construction}[thm]{Construction}
\numberwithin{thm}{section} \numberwithin{equation}{section}

\newcommand{\e}{\epsilon}

\newcommand{\ga}[2]{\begin{gather}\label{#1}#2 \end{gather}}

\newcommand{\surj}{\twoheadrightarrow}

\newcommand{\sD}{{\mathcal D}}
\newcommand{\sE}{{\mathcal E}}
\newcommand{\sF}{{\mathcal F}}
\newcommand{\sG}{{\mathcal G}}
\newcommand{\sH}{{\mathcal H}}

\newcommand{\sL}{{\mathcal L}}

\newcommand{\sN}{{\mathcal N}}
\newcommand{\sO}{{\mathcal O}}
\newcommand{\sP}{{\mathcal P}}

\newcommand{\sR}{{\mathcal R}}
\newcommand{\sS}{{\mathcal S}}

\newcommand{\sU}{{\mathcal U}}
\newcommand{\sV}{{\mathcal V}}

\begin{document}
\title{Remarks on lines
and minimal rational curves}
\author{Ngaiming Mok}
\author{Xiaotao Sun}
\address{Department of Mathematics, The university of Hong Kong, Pokfulam Road, Hong Kong}
\email{nmok@hkucc.hku.hk}
\address{Chinese Academy of Mathematics and Systems Science, Beijing, P. R. of China}
\email{xsun@math.ac.cn}
\date{September 15, 2007}
\thanks{This work is supported by the CERG grant HKU7025/03P of the RGC, Hong Kong.}
\begin{abstract} We determine all of lines in the moduli space $M$
of stable bundles for arbitrary rank and degree. A further
application of minimal rational curves is also given in last
section.
\end{abstract}


\maketitle
\begin{quote}
\end{quote}
\section{Introduction}

Let $C$ be a smooth projective curve of genus $g\ge 2$ and $\sL$
be a line bundle of degree $d$ on $C$. Let $M:=\sS U_C(r,\sL)$ be
the moduli space of stable vector bundles on $C$ of rank $r$ and
with fixed determinant $\sL$, which is a smooth quasi-projective
Fano variety with $Pic(M)=\mathbb Z\cdot \Theta$ and
$-K_M=2(r,d)\Theta$, where $\Theta$ is an ample divisor. In
\cite{Sun}, the second author proved that any rational curve
$\phi:\mathbb P^1\to M$ is defined by a vector bundle $E$ on
$C\times\mathbb P^1$ and gave a formula of its $(-K_M)$-degree in
terms of splitting type of $E$ on the general fiber of
$f:X=C\times \mathbb P^1\to C$. This formula implies immediately
that a rational curve through a general point of $M$ has
$(-K_M)$-degree at least $2r$ and it has degree $2r$ if and only
if it is a Hecke curve. In particular, rational curves of
$(-K_M)$-degree smaller than $2r$, which we call $small$
$rational$ $curves$, must fall in a proper closed subvariety of
$M$. In fact, the formula contains the following information about
points of small rational curves: There exist, for any small
rational curve, a sequence of fixed bundles $F_1,\, F_2,\,
...,\,F_n$ on $C$ such that bundles corresponding points of the
small rational curve are obtained by extensions of $F_1,\, F_2,\,
...,\,F_n$. We should remark here that the bundles $F_1,\, F_2,\,
...,\,F_n$ are independent of points of the small rational curve,
and sometime only depend on the degree of the small rational
curves.

In this paper, we study the rational curves of degree $1$ with
respect to $\Theta$ for arbitrary $r$ and $d$, which we call lines
of $M$. The geometry of $M$ at the case when $(r,d)<r$ is
different from the case when $(r,d)=r$. When $(r,d)<r$, the lines
of $M$ fill up a proper closed subvariety. However, when
$(r,d)=r$, $M$ is generally covered by lines. In Section 2, we
recall firstly two constructions of lines, then, in Theorem
\ref{thm1.7}, we prove that all lines in $M$ are obtained by the
two constructions. In Section 3, we determine the variety
$\bold{Hom}_1(\mathbb P^1,M)$ of degree $1$ morphisms
$\phi:\mathbb P^1\to M$ (Theorem \ref{thm2.1}) and the variety
$\mathscr{L}(M)$ of lines in $M$ (as subvarieties of Chow variety
of $M$) in Corollary \ref{cor2.3}. In Section 4, we present some
partial results on geometry of lines in $M$. In \cite{Sun}, the
proof of main theorem has some implications about the properties
of the bundle $E$ on $C\times \mathbb P^1$ (which defines the
minimal rational curve). In Section 5, we write down first of all
these implications (Theorem \ref{thm4.1}). Then, as an application
of it, we give an alternate proof of some known results (Theorem
\ref{thm4.2}).

\section{The constructions of lines}

Let $C$ be a smooth projective curve of genus $g\ge 2$ and $\sL$ a
line bundle on $C$ of degree $d$. Let $M=\sS\sU_C(r,\sL)^s$ be the
moduli spaces of stable bundles on $C$ of rank $r$, with fixed
determinant $\sL$. It is well-known that $Pic(M)=\mathbb Z\cdot
\Theta$, where $\Theta$ is an ample divisor.

\begin{defn}\label{defn1.1} For any rational curve
$\phi:\mathbb P^1\to M$, its degree is defined to be ${\rm
deg}(\phi^*(\Theta))$. The images $\phi(\mathbb P^1)\subset M$ of
degree $1$ rational curves $\phi:\mathbb P^1\to M$ are called
lines in $M$.
\end{defn}

In this section, we give the constructions of all lines in $M$.
Before stating the first construction, we need the following
lemma, which is a generalization of \cite[Lemma 3.1]{Sun}.

\begin{lem}\label{lem1.2} Let $0\to V_1\to V\to V_2\to 0$ be a nontrivial
extension of vector bundles on $C$. Let $r_i={\rm rk}(V_i)$,
$d_i={\rm deg}(V_i)$ $(i=1,\,2)$, $r={\rm rk}(V)$, $d={\rm
deg}(V)$ be the rank and degree respectively. Then, when
$\,r_1d-d_1r=(r,d)$, $\,V$ is stable if and only if $V_1$ and
$V_2$ are stable.
\end{lem}

\begin{proof} It is clear that $r_i$, $d_i$, $r$, $d$ satisfy
$d_2r-r_2d=(r,d)$, and
$$\mu(V_1)=\mu(V)-\frac{(r,d)}{r_1r},\quad
\mu(V_2)=\mu(V)+\frac{(r,d)}{r_2r}.$$ Writing
$$r_1\frac d{(r,d)} - d_1\frac r{(r,d)} = 1, \quad
d_2\frac r{(r,d)} - r_2\frac d{(r,d)} = 1,$$ we observe that
$(r_i,d_i)=1$ $(i=1,\,2)$.

Assuming that $V$ is stable, we are going to prove the stability
of $V_1$ and $V_2$. For any subbundle $V_1'\subset V_1$ of rank
$r'_1$ and degree $d_1'$, we have, by stability of $V$,
$$r_1'r(\mu(V)-\mu(V_1'))=r_1'd-rd_1'\ge (r,d).$$
Thus $\mu(V_1')\le
\mu(V)-\frac{(r,d)}{r'_1r}=\mu(V_1)+\frac{(r,d)}{r_1r}-\frac{(r,d)}{r'_1r}<\mu(V_1)$,
i.e., $V_1$ is stable. For any subbundle $V_2'\subset V_2$ of rank
$r_2'$, define the subsheaf $V'\subset V$ by $0\to V_1\to V'\to
V_2'\to 0.$ Then $$\mu(V')\le \mu(V)-\frac{(r,d)}{r(r_1+r_2')}$$
and stability of $V_2$ can be seen as follows.
$$\aligned\mu(V_2')&=\mu(V')\frac{r_1+r_2'}{r_2'}-\mu(V_1)\frac{r_1}{r_2'}\\
&\le\mu(V)\frac{r_1+r_2'}{r_2'}-\frac{(r,d)}{r_2'r}-\mu(V_1)\frac{r_1}{r_2'}\\
&=\mu(V)<\mu(V_2).
\endaligned$$

Conversely, assuming that $V_1$ and $V_2$ are stable, we are going
to prove the stability of $V$. For any nontrivial subbundle
$V'\subset V$ of rank $r'$, let $V_2'\subset V_2$ be the image of
$V'$ and $V_1'\subset V_1$ such that
$$0\to V_1'\to V'\to V_2'\to 0$$
is exact. When $V_2'=0$, it is clear that $\mu(V')<\mu(V)$ since
$V_1$ is stable and $\mu(V_1)<\mu(V)$. If $V_1'=0$, then $V_2'$ is
a proper subsheaf of $V_2$ since the extension is nontrivial. Thus
$$\mu(V_2')-\mu(V_2)=-\frac{\text{deg}(V_2^{'*}\otimes
V_2)}{r_2r_2'}<0$$ since $V_2$ is stable. Let $r_2'$, $d_2'$ be
the rank and degree of $V_2'$. Then
$$\aligned\mu(V')&=\mu(V_2')=\mu(V)+\frac{(r,d)}{r_2r}-\frac{\text{deg}(V_2^{'*}\otimes
V_2)}{r_2r_2'}\\&=\mu(V)+\frac{(r,d)}{r_2r_2'r}\left(r_2'-\frac{r}{(r,d)}\text{deg}(V_2^{'*}\otimes
V_2)\right)<\mu(V).\endaligned$$ The last inequality holds because
$r_2'-\frac{r}{(r,d)}\text{deg}(V_2^{'*}\otimes V_2)<r_2$ and is
divisible by $r_2$, thus it must be negative. It is divisible by
$r_2$ since
$$d_2\left(r_2'-\frac{r}{(r,d)}\text{deg}(V_2^{'*}\otimes
V_2)\right)=r_2\left(d_2'-\frac{d}{(r,d)}\text{deg}(V_2^{'*}\otimes
V_2)\right)$$ and $(r_2,d_2)=1$. If $V_1'$, $V_2'$ are nontrivial
of rank $r_1'$, $r'_2$ and degree $d_1'$, $d_2'$, then
$$\aligned\mu(V')&=\mu(V'_1)\frac{r'_1}{r'}+\mu(V_2')\frac{r_2'}{r'}
\le\mu(V_1)\frac{r'_1}{r'}+\mu(V'_2)\frac{r_2'}{r'}\\
&<\mu(V)\frac{r'_1}{r'}+\mu(V_2)\frac{r_2'}{r'}
-\frac{\text{deg}(V_2^{'*}\otimes V_2)}{r_2r'}\\
&=\mu(V)+\frac{(r,d)}{r_2r'r}\left(r_2'-\frac{r}{(r,d)}\text{deg}(V_2^{'*}\otimes
V_2)\right)<\mu(V).\endaligned$$ Thus $V$ is a stable vector
bundle, as desired.
\end{proof}

Now we can describe the first construction of lines. For any given
$r$ and $d$, let $r_1$, $r_2$ be positive integers and $d_1$,
$d_2$ be integers that satisfy the equalities $r_1+r_2=r$,
$d_1+d_2=d$ and
$$r_1\frac{d}{(r,d)}-d_1\frac{r}{(r,d)}=1,\quad d_2\frac{r}{(r,d)}-r_2\frac{d}{(r,d)}=1.$$
Let $\sU_C(r_1,d_1)$ (resp. $\sU_C(r_2,d_2)$) be the moduli space
of stable vector bundles of rank $r_1$ (resp. $r_2$) and degree
$d_1$ (resp. $d_2$). Then, since $(r_1,d_1)=1$ and $(r_2,d_2)=1$,
they are smooth projective varieties and there are universal
vector bundles $\sV_1$, $\sV_2$ on $C\times \sU_C(r_1,d_1)$ and
$C\times\sU_C(r_2,d_2)$ respectively. Consider the morphism
$$\sU_C(r_1,d_1)\times\sU_C(r_2,d_2)\xrightarrow{{\rm det}(\bullet)\times
{\rm det}(\bullet)} J^{d_1}_C\times
J^{d_2}_C\xrightarrow{(\bullet)\otimes (\bullet)}J^d_C$$ and let
$\sR(r_1,d_1)$ be its fiber at $[\sL]\in J^d_C$. We still use
$\sV_1$, $\sV_2$ to denote the pullback on $C\times \sR(r_1,d_1)$
by the projection $C\times\sR(r_1,d_1)\to C\times\sU_C(r_i,d_i)$
($i=1,2$) respectively. Let $p:C\times\sR(r_1,d_1)\to
\sR(r_1,d_1)$ and $\sG=R^1p_*(\sV_2^{\vee}\otimes\sV_1)$. Then,
since ${\rm Hom}(V_2,V_1)=0$, $\sG$ is a vector bundle of rank
$r_1r_2(g-1)+(r,d)$. Let $q: P(r_1,d_1)=\mathbb
P(\sG)\to\sR(r_1,d_1)$ be the projective bundle parametrzing
$1$-dimensional subspaces of $\sG_t$ ($t\in\sR(r_1,d_1)$) and
$f:C\times P(r_1,d_1)\to C$, $\pi:C\times P(r_1,d_1)\to
P(r_1,d_1)$ be the projections. Then there is a universal
extension \ga{1.1}{0\to (id\times
q)^*\sV_1\otimes\pi^*\sO_{P(r_1,d_1)}(1)\to\sE\to (id\times
q)^*\sV_2 \to 0} on $C\times P(r_1,d_1)$ such that for any
$x=([V_1], [V_2], [e])\in P(r_1,d_1)$, where $[V_i]\in
\sU_C(r_i,d_i)$ with ${\rm det}(V_1)\otimes {\rm det}(V_2)=\sL$
and $[e]\subset {\rm H}^1(C,V_2^{\vee}\otimes V_1)$ being a line
through the origin, the bundle $\sE|_{C\times\{x\}}$ is the
isomorphic class of vector bundles $E$ given by extensions
$$0\to V_1\to E\to V_2\to 0 $$
that defined by vectors on the line $[e]\subset {\rm
H}^1(C,V_2^{\vee}\otimes V_1)$.

To see the existence of the universal extension, recall
\cite[Lemma 2.4]{Ra}: For two families $(E_s)_{s\in S}$,
$(F_t)_{t\in T}$ of bundles on $C\times S$, $C\times T$, there
exists a universal extension if (1) ${\rm dim}\,{\rm H}^1(C,\sH
om(F_t,E_s))$ is independent of $(s,t)\in S\times T$, (2) ${\rm
H}^i(S\times T, p_{S\times T*}\sH om(F,E)\otimes V^*)=0$
($i=1,2$), where $V$ is the vector bundle on $S\times T$ with
fibers ${\rm H}^1(C,\sH om(F_t,E_s))$ at $(s,t)\in S\times T$. In
our case, $E=\sV_1$, $F=\sV_2$, and the above conditions are
satisfied since ${\rm Hom}(V_2,V_1)=0$ for any $[V_i]\in
\sU_C(r_i,d_i)$ ($i=1,2$). By Lemma \ref{lem1.2}, the universal
extension $$0\to (id\times
q)^*\sV_1\otimes\pi^*\sO_{P(r_1,d_1)}(1)\to\sE\to (id\times
q)^*\sV_2 \to 0$$ on $C\times P(r_1,d_1)$ defines a morphism
\ga{1.2}{\Phi:P(r_1,d_1)\to \sS\sU_C(r,\sL)^s=M.}

\begin{construction} The images (under $\Phi$) of lines in the fibres of
$$q: P(r_1,d_1)=\mathbb P(\sG)\to\sR(r_1,d_1)$$
are lines of $M$.
\end{construction}

\begin{lem}\label{lem1.4} On each fiber $P(r_1,d_1)_{\xi}:=q^{-1}(\xi)$ at $\xi\in
\sR(r_1,d_1)$,
$$\Phi_{\xi}:=\Phi|_{P(r_1,d_1)_{\xi}}:P(r_1,d_1)_{\xi}\to M$$
is the normalization of its image. In particular, the rational
curves constructed in  \textbf{Construction 2.3} are lines in $M$.
\end{lem}

\begin{proof} Write $\sP=P(r_1,d_1)=\mathbb P(\sG)$, $\sR=\sR(r_1,d_1)$,
$\sU_i=\sU_C(r_i,d_i)$ ($i=1,2$), recall
$\sG=R^1p_*(\sV_2^{\vee}\otimes\sV_1)$.  Let
$\omega_C=\sO_C(\sum^{2g-2}_{i=1}y_i)$, $\omega_{\sU_1}$ and
$\omega_{\sU_2}$ be the canonical line bundles of $C$, $\sU_1$ and
$\sU_2$. It is not difficult, using \eqref{1.1}, to compute
$\Phi^*(\omega^{-1}_M)=$
$$q^*\left(\omega_{\sU_1}^{-1}\otimes\omega_{\sU_2}^{-1}\otimes{\rm det}(\sG)^{\otimes 2}\otimes
\bigotimes^{2g-2}_{i=1}{\rm
det}(\sV_1^{\vee}\otimes\sV_2)_{y_i}\right)\otimes
\sO_{\sP}(2(r,d)).$$ Thus, for any $\xi\in\sR$,
$\Phi_{\xi}^*(\Theta)=\sO_{\sP_{\xi}}(1)$ and $\Phi_{\xi}$ is the
normalization of $\Phi_{\xi}(\sP_{\xi})$. In particular, for any
line $\ell\subset \sP_{\xi}$, $\Phi_{\xi}(\ell)\subset M$ is a
line and $\ell$ is the normalization of $\Phi_{\xi}(\ell)$.
\end{proof}

Now we recall the construction of Hecke curves which are also
lines in $M$ when $(r,d)=r$. Let $\sU_C(r,d-1)$ be the moduli
space of stable bundles of rank $r$ and degree $d-1$. Let
$\mathfrak{O}\subset\sU_C(r,d-1)$ be the open set of
$(1,0)$-semistable bundles in the following sense

\begin{defn}\label{defn1.5} A vector bundle $V$ on $C$ is called
$(k,\ell)$-semistable (resp. $(k,\ell)$-stable) if for any proper
subbundle $W\subset V$, we have
$$ \frac{{\rm deg}(W)+k}{{\rm rk}(W)}\le\,(\text{resp. $<$})
\frac{{\rm deg}(V)+k-\ell}{{\rm rk}(V)}.$$
\end{defn}

Let $C\times\mathfrak{O}\xrightarrow\psi J^d(C)$ be defined as
$\psi(x, V)=\sO_C(x)\otimes{\rm det}(V)$ and let
$\mathscr{R}_C:=\psi^{-1}(\sL)\subset C\times \mathfrak{O}$. There
is a fibration $\mathscr{R}_C\to C$ with fibres
$\mathfrak{O}\cap\sS\sU_C(r,\sL(-x))$ at $ x\in C$. Let
$\mathscr{V}$ be the universal bundle on $\mathscr{R}_C$, let
$p:\mathbb P(\mathscr{V}^{\vee})\to \mathscr{R}_C$ be the
projective bundle and
$$p^*(\mathscr{V}^{\vee})\to\sO_{\mathbb P(\mathscr{V}^{\vee})}(1)\to 0$$
be the universal quotient. Let $C\times \mathbb
P(\mathscr{V}^{\vee})\xrightarrow\pi \mathbb
P(\mathscr{V}^{\vee})$ be the projection and $\Gamma\subset
C\times \mathbb P(\mathscr{V}^{\vee})$ be the graph of $\mathbb
P(\mathscr{V}^{\vee})\xrightarrow p\mathscr{R}_C\to C$. We have
$$0\to\mathscr{E}^{\vee}\to\pi^*p^*(\mathscr{V}^{\vee})\to \sO_{\Gamma}\otimes\pi^*\sO_{\mathbb
P(\mathscr{V}^{\vee})}(1)\to 0$$ where $\mathscr{E}^{\vee}$ is
defined to the kernel of the surjection. Taking duals, we have
\ga{1.3} {0\to
\pi^*p^*\mathscr{V}\to\mathscr{E}\to\sO_{\Gamma}(\Gamma)\otimes\pi^*\sO_{\mathbb
P(\mathscr{V}^{\vee})}(-1)\to 0,} which, at any point $\xi=(x, V,
V_x^{\vee}\surj \Lambda)\in \mathbb P(\mathscr{V}^{\vee})$, gives
exact sequence $$0\to V\xrightarrow\iota\mathscr{E}_{\xi}\to
\,\sO_x\to 0$$ on $C$ such that ${\rm
ker}(\iota_x)=\Lambda^{\vee}\subset V_x$. That $V$ being
$(1,0)$-semistable (resp. stable) implies semistability (resp.
stability) of $\mathscr{E}_{\xi}$ defines
$$\Psi:\mathbb P(\mathscr{V}^{\vee})\to \sS\sU_C(r,\sL)\supseteq\sS\sU_C(r,\sL)^s=M.$$
Let $\mathscr{P}^0:=\Psi^{-1}(M)\subset\mathbb
P(\mathscr{V}^{\vee})$,
$\mathscr{R}_C^0:=p(\mathscr{P}^0)\subset\mathscr{R}_C$ and
\ga{1.4}{p: \mathscr{P}^0\to \mathscr{R}_C^0\,,\quad \Psi:
\mathscr{P}^0\to M\,.}

\begin{construction} The images (under $\Psi$) of lines in the fibres of
$$p: \mathscr{P}^0\to \mathscr{R}_C^0$$
are the so called Hecke curves in $M$, which are lines if and only
if $(r,d)=r$ by \cite [Theorem 1]{Sun}.
\end{construction}

\begin{thm}\label{thm1.7} (i) When $(r,d)\neq r$, all lines in $M$ are obtained by performing
\textbf{Construction 2.3} for all pairs $\{r_1,d_1\}$ satisfying
$$0<r_1<r,\quad r_1d-d_1r=(r,d).$$
(ii) When $(r,d)=r$, perform \textbf{Construction 2.3} as in (i)
and the \textbf{Construction 2.6}, we obtain all lines in $M$.
\end{thm}

\begin{proof} In \cite{Sun} it was shown that any rational curve
$\phi:\mathbb P^1\to M$ is defined by a vector bundle $E$ on
$X=C\times\mathbb P^1$ and also have proven a degree formula. To
recall it, let $f: X\to C$ and $\pi: X\to\mathbb P^1$ be the
projections.  On a general fiber $f^{-1}(\xi)=X_{\xi}$, $E$ has
the form
$$E|_{X_{\xi}}=\bigoplus_{i=1}^n\sO_{X_{\xi}}(\alpha_i)^{\oplus r_i},\quad
\alpha_1>\,\,\cdots\,\,>\alpha_n.$$  The $\alpha=(\alpha_1^{\oplus
r_1},\,...,\,\alpha_n^{\oplus r_n})$ is called the generic
splitting type of $E$. Tensoring $E$ by $\pi^*\sO_{\mathbb P^1}
(-\alpha_n)$, we can assume without loss of generality that
$\alpha_n=0$. Any such $E$ admits a relative Harder-Narasimhan
filtration
$$0=E_0\subset E_1\subset \,\cdots\,\subset E_n=E$$
in which the quotient sheaves $F_i=E_i/E_{i-1}$ are torsion-free
with generic splitting type $(\alpha_i^{\oplus r_i})$
respectively. Let $F_i'=F_i\otimes \pi^*\sO_{\mathbb
P^1}(-\alpha_i)$ ($i=1,\,...,n$), thus they have generic splitting
type $(0^{\oplus r_i})$ respectively. Without risk of confusion,
we denote the degree of $F_i$ (resp. $E_i$) on the general fiber
of $\pi$ by $\text{deg}(F_i)$ (resp. $\text{deg}(E_i)$).
Accordingly, $\mu(E_i)$ (resp. $\mu(E)$) denotes the slope of the
restriction of $E_i$ (resp. $E$) to the general fiber of $\pi$
respectively. Let $\text{rk}(E_i)$ denote the rank of $E_i$. Then
we have the formula (See the formula (2.2) of \cite{Sun})
$$\text{deg}(\phi^*(\Theta))=\frac{r}{(r,d)}\left(\sum_{i=1}^nc_2(F'_i)+\sum_{i=1}^{n-1}
(\mu(E)-\mu(E_i))(\alpha_i-\alpha_{i+1})\text{rk}(E_i) \right).$$

When $(r,d)\neq r$,  we have $c_2(F'_i)=0$ and $n=2$. Thus there
are bundles $V_1$, $V_2$ of rank $r_1$, $r_2$ and degree $d_1$,
$d_2$ on $C$ such that
$$0\to f^*V_1\otimes\pi^*\sO_{\mathbb P^1}(1)\to E\to f^*V_2\to 0$$
where $r_1$, $r_2$, $d_1$, $d_2$ satisfy $r_1+r_2=r$, $d_1+d_2=d$
and
$$r_1\frac{d}{(r,d)}-d_1\frac{r}{(r,d)}=1,\quad d_2\frac{r}{(r,d)}-r_2\frac{d}{(r,d)}=1.$$
By Lemma \ref{lem1.2}, $V_1$ and $V_2$ must be stable and
$det(V_1)\otimes det(V_2)=\sL$. Thus $\phi$ factors through
$\mathbb
P^1\xrightarrow{\sigma}\sP_{\xi}\xrightarrow{\Phi_{\xi}}M$, where
$\xi=(V_1,V_2)\in\sR$ and $\sigma^*\sO_{\sP_{\xi}}(1)=\sO_{\mathbb
P^1}(1)$ (so that $\sigma$ is an embedding and $\sigma(\mathbb
P^1)$ is a line of $\sP_{\xi}$). This proves (i).

When $(r,d)=r$, we have either $c_2(F'_i)=0$ and $n=2$ or
$c_2(E)=1$ and $n=1$. Thus the line is either obtaining by
\textbf{Construction 2.3} or defined by a vector bundle $E$ on
$X=C\times\mathbb P^1$ satisfying  $$0\to f^*V\to
E\to\sO_{\{p\}\times \mathbb P^1}(-1)\to 0$$ where
$f:X=C\times\mathbb P^1\to C$ and $\pi:X=C\times\mathbb
P^1\to\mathbb P^1$ are projections, $V$ is a vector bundle on $C$.
The stability of $E_t=E|_{C\times\{t\}}$ ($\forall\,\, t\in\mathbb
P^1$) implies immediately that $V$ is $(1,0)$-semistable. Thus, in
this case, the line is obtained by \textbf{Construction 2.6}.
\end{proof}

\section{The variety of lines}

By the variety of lines, we mean the quotient
$\bold{Hom}_1(\mathbb P^1,M)/{\rm Aut}(\mathbb P^1)$ which can be
defined by means of the Chow variety. To determine
$\bold{Hom}_1(\mathbb P^1,M)$, recall from \textbf{Construction
2.3} and \textbf{Construction 2.6}, we have
$$q: P(r_1,d_1)=\mathbb P(\sG)\to\sR(r_1,d_1),\quad
p: \mathscr{P}^0\to \mathscr{R}_C^0.$$ Let $\mathbb
P^1_{\sR(r_1,d_1)}=\mathbb P^1\times \sR(r_1,d_1)$ and
$\bold{Hom}_1\left(\mathbb P^1_{\sR(r_1,d_1)}\,,\,\,\mathbb
P(\sG)/\sR(r_1,d_1)\right)$ be the scheme such that for any scheme
$T$ over $\sR(r_1,d_1)$
$$\bold{Hom}_1\left(\mathbb
P^1_{\sR(r_1,d_1)}\,,\,\,\mathbb P(\sG)/\sR(r_1,d_1)\right)(T)$$
is the set of $T$-morphisms $\mathbb
P^1_{\sR(r_1,d_1)}\times_{\sR(r_1,d_1)}T\to
P(\sG)\times_{\sR(r_1,d_1)}T$ of degree $1$ with respect to
$\sO_{P(\sG)}(1)$. It is the variety of degree $1$ maps
$$\mathbb P^1\to P(r_1,d_1)=\mathbb P(\sG)$$ with images in the fibers of
$q:P(r_1,d_1)=\mathbb P(\sG)\to\sR(r_1,d_1)$. Similarly, recall
that $p: \mathscr{P}^0\to \mathscr{R}_C^0$ is an open set of the
projective bundle $p:\mathbb P(\mathscr{V}^{\vee})\to
\mathscr{R}_C$, we can define the variety
$$\bold{Hom}^r_1(\mathbb P^1,\mathscr{P}^0):=\bold{Hom}_1\left(\mathbb
P^1_{\mathscr{R}_C^0}\,,\,\,\mathscr{P}^0/\mathscr{R}_C^0\right)$$
of degree $1$ maps $\mathbb P^1\to\mathscr{P}^0$ with images in
the fibers of $p: \mathscr{P}^0\to \mathscr{R}_C^0$ (we use
$\bold{Hom}^r$ to denote relative maps). Let
$$\bold{Hom}^r_1(\mathbb P^1,\mathbb P):=\bigsqcup_{\{r_1,d_1\}} \bold{Hom}_1\left(\mathbb
P^1_{\sR(r_1,d_1)}\,,\,\,\mathbb P(\sG)/\sR(r_1,d_1)\right) $$ be
the disjoint union, where $\{r_1,d_1\}$ runs through the pairs
satisfying:
$$0<r_1<r,\quad r_1d-d_1r=(r,d).$$

\begin{thm}\label{thm2.1} Let $\bold{Hom}_1(\mathbb P^1,M)$ be the
variety of degree $1$ morphisms $\mathbb P^1\to M$ (respect to
$\Theta$). Then
$$\bold{Hom}_1(\mathbb P^1,M)\cong\begin{cases}\bold{Hom}^r_1(\mathbb P^1,\mathbb P)  & \text{if $(r,d)\neq r$}\\
\bold{Hom}^r_1(\mathbb P^1,\mathbb P)\bigsqcup
\bold{Hom}^r_1(\mathbb P^1,\mathscr{P}^0)&\text{if
$(r,d)=r$}.\end{cases}$$
\end{thm}

\begin{proof} By sending a $T$-morphism $$\mathbb P^1\times T\cong\mathbb
P^1_{\sR(r_1,d_1)}\times_{\sR(r_1,d_1)}T\xrightarrow{\varphi_T}
P(\sG)\times_{\sR(r_1,d_1)}T$$ to a $T$-morphism
$$\mathbb P^1\times T\xrightarrow{\varphi_T} P(\sG)\times_{\sR(r_1,d_1)}T\to
P(\sG)\times T\xrightarrow{\Phi\times id_T}M\times T,$$ we have
the canonical morphism
$$\bold{Hom}^r_1(\mathbb P^1,\mathbb
P)\to \bold{Hom}_1(\mathbb P^1,M)$$ which is surjective when
$(r,d)\neq r$ by Theorem \ref{thm1.7}. To show it is also
injective, let $\xi_1,\,\xi_2\in\bold{Hom}^r_1(\mathbb P^1,\mathbb
P)$ defined by the exact sequences
$$0\to f^*V_1\otimes\pi^*\sO_{\mathbb P^1}(1)\to\sE_1\to f^*V_2\to
0,\,$$  $$ 0\to f^*W_1\otimes\pi^*\sO_{\mathbb P^1}(1)\to\sE_2\to
f^*W_2\to 0$$ on $C\times\mathbb P^1$, where $f:X=C\times\mathbb
P^1\to C$ and $\pi:X=C\times\mathbb P^1\to\mathbb P^1$ are the
projections. If $\xi_1$, $\xi_2$ have the same image in
$\bold{Hom}_1(\mathbb P^1,M)$, then there is a line bundle $\sN$
on $\mathbb P^1$ such that $\sE_1\cong \sE_2\otimes\pi^*\sN$. If
$\text{deg}(\sN)\le 0$, then
$\text{Hom}(f^*V_1\otimes\pi^*\sO_{\mathbb
P^1}(1),f^*W_2\otimes\pi^*\sN)=0$ and $\sE_1\cong
\sE_2\otimes\pi^*\sN$ induces $f^*V_1\otimes\pi^*\sO_{\mathbb
P^1}(1)\hookrightarrow f^*W_1\otimes\pi^*\sO_{\mathbb
P^1}(1)\otimes\pi^*\sN$, which implies $V_1\hookrightarrow
W_1\otimes\text{H}^0(\sN)$. Thus $\sN=\sO_{\mathbb P^1}$,
$V_1\cong W_1$ and $V_2\cong W_2$, which implies $\xi_1=\xi_2$. If
$\text{deg}(\sN)\ge 0$, using
$\sE_2\cong\sE_1\otimes\pi^*\sN^{-1}$, we get $\xi_1=\xi_2$ by the
same arguments. Thus $\bold{Hom}^r_1(\mathbb P^1,\mathbb P)\to
\bold{Hom}_1(\mathbb P^1,M)$ is bijective when $(r,d)\neq r$.

Similarly, when $(r,d)=r$,  we have a surjective morphism
$$\bold{Hom}^r_1(\mathbb P^1,\mathbb P)\bigsqcup
\bold{Hom}^r_1(\mathbb P^1,\mathscr{P}^0)\to\bold{Hom}_1(\mathbb
P^1,M)$$ by Theorem \ref{thm1.7}. To see the injectivity, we only
need to consider $\xi_1,\,\xi_2\in\bold{Hom}^r_1(\mathbb
P^1,\mathscr{P}^0)$ defined by the following two exact sequences
on $C\times \mathbb P^1$
$$0\to f^*V\to \sE_1\to\sO_{\{x_1\}\times \mathbb P^1}(-1)\to 0,$$
$$0\to f^*W\to\sE_2\to\sO_{\{x_2\}\times \mathbb P^1}(-1)\to 0$$ where $V$, $W$
are stable vector bundles on $C$ of rank $r$ and degree $d-1$,
$x_1,\,x_2\in C$ are two points. If $\xi_1$, $\xi_2$ have the same
image in $\bold{Hom}_1(\mathbb P^1,M)$, then there is a line
bundle $\sN$ on $\mathbb P^1$ such that
$$\sE_1\cong \sE_2\otimes\pi^*\sN,\quad x_1=x_2.$$
If ${\rm deg}(\sN)\le  0$, then ${\rm Hom}(f^*V,\sO_{\{x_2\}\times
\mathbb P^1}(-1)\otimes \pi^*\sN)=0$. The isomorphism $\sE_1\to
\sE_2\otimes\sN$ induces an injection $f^*V\hookrightarrow
f^*W\otimes\pi^*\sN$, which implies that $\sN=\sO_{\mathbb P^1}$
and $V\cong W$, thus $\xi_1=\xi_2$. If $\text{deg}(\sN)\ge 0$,
using $\sE_2\cong\sE_1\otimes\pi^*\sN^{-1}$, we have $\xi_1=\xi_2$
by the same arguments.

To show the isomorphism, it is enough to show that
$\bold{Hom}_1(\mathbb P^1,M)$ is smooth. To see the smoothness of
$\bold{Hom}_1(\mathbb P^1,M)$, let $\varphi:\mathbb P^1\to M$ be a
point of $\bold{Hom}_1(\mathbb P^1,M)$, which, by Lemma 2.1 of
\cite{Sun}, is defined by a vector bundle $E$ on $C\times\mathbb
P^1$ such that $\varphi^*T_M=R^1\pi_*Ad(E)$. Then $E$ must satisfy
either $0\to f^*V_1\otimes\pi^*\sO_{\mathbb P^1}(1)\to E\to
f^*V_2\to 0\,\,$ or
$$0\to f^*V\to E\to\sO_{\{x\}\times \mathbb P^1}(-1)\to 0.$$ Using these
exact sequences, we can show
$${\rm H}^1(\varphi^*T_M)={\rm H}^1(R^1\pi_*Ad(E))=0.$$ Thus
$\bold{Hom}_1(\mathbb P^1,M)$ is smooth.
\end{proof}

By Theorem 3.21 of \cite{Ko}, there is a semi-normal variety ${\rm
Chow}_{1,1}(M)$ parametrizing effective cycles of dimension $1$
and degree $1$ (respect to $\Theta$) with a universal cycle ${\rm
Univ}_{1,1}(M)\to {\rm Chow}_{1,1}(M).$ Since
$\bold{Hom}_1(\mathbb P^1,M)$ is smooth, there is an ${\rm
Aut}(\mathbb P^1)$-invariant morphism
$$\bold{Hom}_1(\mathbb P^1,M)\to {\rm Chow}_{1,1}(M).$$ Let
$\mathscr{L}(M)\subset{\rm Chow}_{1,1}(M)$ be the image, which is
precisely the locus of ${\rm Chow}_{1,1}(M)$ parametrizing the
cycles with rational components. Then, by Proposition 2.2 of
\cite{Ko}, $\mathscr{L}(M)\subset{\rm Chow}_{1,1}(M)$ is a closed
subset.

\begin{defn}\label{defn2.2} The closed subset $\mathscr{L}(M)\subset{\rm
Chow}_{1,1}(M)$ with the reduced scheme structure is called the
\textbf{variety of lines in $M$}. The induced universal cycle
$\mathbb{L}\subset M\times\mathscr{L}(M)$ defined by
$$\mathbb{L}:={\rm Univ}_{1,1}(M)\times_{{\rm
Chow}_{1,1}(M)}\mathscr{L}(M)\to\mathscr{L}(M)$$ is called the
\textbf{universal line in $M$}.
\end{defn}

Let $G(r_1,d_1)\to\sR(r_1,d_1)$ (resp. $\mathscr{H}\to
\mathscr{R}_C^0$) be the relative Grassmannian  bundles of lines
in $P(r_1,d_1)$ (resp. $\mathscr{P}^0$), and let
$$\xymatrix{\mathfrak{L}(r_1,d_1)\ar[dr]\ar@{^{(}->}[r]
& P(r_1,d_1)\times_{\sR(r_1,d_1)}G(r_1,d_1) \ar[d]  \\
&    G(r_1,d_1)} \,,\quad\xymatrix{ \mathfrak{L}(h) \ar[dr]
\ar@{^{(}->}[r]& \mathscr{P}^0\times_{\mathscr{R}_C^0}\mathscr{H}\ar[d]  \\
& \mathscr{H}}$$ be the universal lines. Recall the morphisms
\ga{2.1}{\Phi: P(r_1,d_1)\to M\,,\quad \Psi:\mathscr{P}^0\to M} in
\eqref{1.2} and \eqref{1.4}, which induce \ga{2.2}{\xymatrix{
\mathfrak{L}(r_1,d_1)\ar[dr]\ar[r]^{\Phi\times {\rm id}}& M\times G(r_1,d_1)\ar[d]  \\
&  G(r_1,d_1)  }\,,\quad \xymatrix{\mathfrak{L}(h) \ar[dr]
\ar[r]^{\Psi\times{\rm id}}& M\times\mathscr{H}\ar[d]  \\
& \mathscr{H}}} Then the families ${\rm Im}(\Phi\times{\rm
id})\subset M\times G(r_1,d_1)$ and ${\rm Im}(\Psi\times{\rm
id})\subset M\times\mathscr{H}$ of lines define the morphisms
\ga{2.3}{G(r_1,d_1)\xrightarrow{\Upsilon_{r_1,d_1}}
\mathscr{L}(M)\quad \text{and\,\,
$\mathscr{H}\xrightarrow\theta\mathscr{L}(M)$\,\, if $(r,d)=r$}.}

Let $\mathscr{L}(M)_{r_1,d_1}:={\rm Im}(\Upsilon_{r_1,d_1})$,
$\mathscr{H}_{\theta}:={\rm Im}(\theta)$, and let
$$G(M):=\bigsqcup _{\{r_1,d_1\}}G(r_1,d_1)\,,\quad \mathscr{S}(M):=
\bigsqcup_{\{r_1,d_1\}}\mathscr{L}(M)_{r_1,d_1}$$ be the disjoint
unions of varieties, where $\{r_1,d_1\}$ runs through the pairs
$\{r_1,d_1\}$ satisfying: $0<r_1<r,\quad r_1d-d_1r=(r,d)$.

\begin{cor}\label{cor2.3}$G(r_1,d_1)\xrightarrow{\Upsilon_{r_1,d_1}}\mathscr{L}(M)_{r_1,d_1}$
are the normalizations and $\theta$ induces
$\mathscr{H}\cong\mathscr{H}_{\theta}\subset\mathscr{L}(M)$ when
$(r,d)=r$. Moreover,
$$\mathscr{L}(M)= \begin{cases} \mathscr{S}(M)  & \text{if $(r,d)\neq r$}\\
\mathscr{S}(M)\sqcup \mathscr{H}_{\theta}  &\text{if
$(r,d)=r$}\end{cases}$$ and $G(M)\to\mathscr{L}(M)$ is an
injective morphism.
\end{cor}

\begin{proof} $\mathscr{H}\cong\mathscr{H}_{\theta}$
follows from the study of Hecke cycles in \cite{NR}. Since all
$G(r_1,d_1)$ are smooth projective varieties, to show the other
statements, it is enough to show that $G(M)\to \mathscr{L}(M)$
(resp. $G(M)\sqcup\mathscr{H}\to\mathscr{L}(M)$) is bijective if
$(r,d)\neq r$ (resp. $(r,d)=r$). Theorem \ref{thm1.7} implies
surjectivity. The same arguments in the proof Theorem \ref{thm2.1}
imply injectivity.
\end{proof}

\section{The geometry of lines}

The morphism $\Psi:\mathscr{P}^0\to M$ was well studied in
\cite{NR} for arbitrary rank. In particular, for any
$\xi\in\mathscr{R}_C^0$, the morphism
$$\Psi_{\xi}:=\Psi|_{\mathscr{P}^0_{\xi}}: \mathscr{P}^0_{\xi}=p^{-1}(\xi)\to M$$
is a closed embedding. In this section, we study the morphism
$$\Phi: \sP:=P(r_1,d_1)\to M$$
for arbitrary rank. In general, we are not able to show that
$$\Phi_{\xi}:=\Phi|_{\sP_{\xi}}:\sP_{\xi}=q^{-1}(\xi)\to M$$
is a closed embedding for each $\xi\in\sR:=\sR(r_1,d_1)$.
Consequently, we are not able to show that every line in $M$ is
smooth for arbitrary rank case (it is true in rank two case).
However, we will show that $\Phi_{\xi}$ is a closed embedding for
$\xi\in \sR\setminus\sD$, where
$$\sD=\{\xi=(V_1,V_2)\in \sR\,|\,{\rm Hom}(V_1,V_2)\neq 0\}.$$
It can be realized as the degeneration locus of a morphism between
two vector bundles on $\sR$. Thus, if $\sD\neq\emptyset$, it has
$${\rm Codim}(\sD)\le r_1r_2(g-1)+1-(r,d)$$
and $\sD$ is Cohen-Macaulay if the equality holds. To prove a
lower bound of the codimension, we start it with an elementary
lemma.

\begin{lem}\label{lem3.1} Let $V_1$, $V_2$, $V$ be stable vector bundles on
$C$ of rank $r_1$, $r_2$, $r=r_1+r_2$ and degree $d_1$, $d_2$,
$d=d_1+d_2$. Then, when $\,r_1d-d_1r=(r,d)$, we have
\begin{itemize}\item[(1)] Any nontrivial morphism $V_1\to V$ must be an
injective morphism of bundles, and any nontrivial morphism $V\to
V_2$ must be surjective.\item[(2)] For any nontrivial morphism
$f:V_1\to V_2$, if $\mu(f(V_1))\neq \mu(V)$, then it must be
injective when $r_1\le r_2$ and surjective when $r_1>r_2$. If
$\mu(f(V_1))=\mu(V)$ and $(r,d)\neq r$, then $f(V_1)$ is
semistable and $V_2/f(V_1)$ is torsion-free.
\end{itemize}
\end{lem}

\begin{proof} Let $V_1'\subset V$ be the image of $V_1\to V$ with ${\rm rk}(V_1')=r_1'$, ${\rm deg}(V_1')=d_1'$. Then
$$\frac{(r,d)}{r_1r}=\mu(V)-\mu(V_1')+\mu(V_1')-\mu(V_1)>\mu(V)-\mu(V_1')=\frac{r_1'd-rd_1'}{r'_1r}>0$$
if $r_1'\neq r_1$, which is impossible since $r_1'd-rd_1'\ge
(r,d)$. It also shows that $V_1'$ must be a subbundle of $V$. The
surjectivity of any nontrivial morphism $V\to V_2$ can be proved
similarly. To prove (2), let $f(V_1)$ be of rank $r_1'$ and degree
$d_1'$, then
\ga{3.1}{\frac{(r,d)}{r_1r_2}=\mu(V_2)-\mu(f(V_1))+\mu(f(V_1))-\mu(V_1)\\=\frac{r_1'd_2-r_2d_1'}{r'_1r_2}+
\frac{r_1d_1'-r_1'd_1}{r'_1r_1}.\notag} When $r_1\le r_2$, if
$V_1\to V_2$ is not injective, then both $\text{deg$(V_2\otimes
f(V_1)^*)$}= r_1'd_2-r_2d_1'$ and $\text{deg$(f(V_1)\otimes
V_1^*)$}= r_1d_1'-r_1'd_1$ are positive. Their difference
$$(r_1'd_2-r_2d_1')-(r_1d_1'-r_1'd_1)=r_1'd-d'_1r=r_1'r(\mu(V)-\mu(f(V_1)))\neq 0$$ is a
nonzero integer divisible by $(r,d)$, thus one of them is bigger
than $(r,d)$, which contradicts the above equality \eqref{3.1}.
When $r_1>r_2$, then $\text{deg$(f(V_1)\otimes V_1^*)$}>0$ and
$\text{deg$(V_2\otimes f(V_1)^*)$}\ge 0$. The same argument shows
that $\text{deg$(V_2\otimes
f(V_1)^*)$}=r_2r_1'(\mu(V_2)-\mu(f(V_1)))$ must be zero. Thus
$f(V_1)=V_2$ by the stability of $V_2$.

When $\mu(f(V_1))=\mu(V)$, we show first of all that $V_2/f(V_1)$
is torsion-free. Let $f(V_1)\subset W\subset V_2$ such that
$V_2/W$ is torsion-free, ${\rm rk}(W)=r_1'$, ${\rm deg}(W)=\tilde
d'_1$. Then $(r_1\tilde d_1'-r_1'd_1) -(r_1'd_2-r_2\tilde
d_1')=r(\tilde d_1'-d_1')$ and
$$\aligned\frac{(r,d)}{r_1r_2}&=\mu(V_2)-\mu(W)+\mu(W)-\mu(V_1)\\
&=\frac{r_1'd_2-r_2\tilde d_1'}{r_1'r_2}+\frac{r_1\tilde
d_1'-r_1'd_1}{r'_1r_1}\ge\frac{r(\tilde
d_1'-d_1')}{r_1r_1'}.\endaligned$$ Thus, if $\tilde d_1'-d_1'>0$,
we get $r\le (r,d)$, which contradicts the assumption $r\neq
(r,d)$. To see that $f(V_1)$ is semistable, let $V_0\subset
f(V_1)$ be a proper subbundle of rank $r_0$ and degree $d_0$. If
$\mu(V_0)>\mu(f(V_1))=\mu(V)$, then $\mu(V_1)<\mu(V_0)<\mu(V_2),$
which is impossible by the above arguments. Thus $f(V_1)$ is
semistable.
\end{proof}

\begin{lem}\label{lem3.2} Let $\sD=\{(V_1, V_2)\in\sU_C(r_1,d_1)\times\sU_C(r_2,d_2)|{\rm Hom}(V_1,V_2)\neq
0\}$ and $\sR:=\sR(r_1,d_1)$. Then, when ${\rm
min}\{r_1,r_2\}>\frac{r}{(r,d)}$, we have \ga{3.2}{{\rm
codim}(\sD\cap\sR)\ge\frac{r}{(r,d)}(r-\frac{r}{(r,d)})(g-1)-1,}
and when ${\rm min}\{r_1,r_2\}\le\frac{r}{(r,d)}$, we have
\ga{3.3}{{\rm codim}(\sD\cap\sR)\ge r_1r_2(g-1)+1-(r,d)} The same
inequalities also hold for the codimension of $\sD$.
\end{lem}

\begin{proof} Since taking dual of vector bundles induces an
isomorphism between moduli spaces, we can assume $r_1\ge r_2$
without loss of generality. Let $h:\sH\to\sD$ be the total space
of morphisms $V_1\to V_2$, let $\sH_1\subset\sH$ be the union of
irreducible components whose general points are not surjective
morphisms $V_1\to V_2$, and $\sH_2:=\sH\setminus\sH_1$. Then there
is an open dense subset $\sH_2^0\subset\sH_2$ and an exact
sequence $0\to\sV'\to\sV_1\to\sV_2\to 0$ on $C\times\sH_2^0$,
where $\sV'$ is a flat family of vector bundles of rank $r_1-r_2$
and degree $d_1-d_2$ any subbundle of which has slope less than
$d_1/r_1$ (so that the set of such bundles is bounded). Let
$Q\subset {\rm Quot}(\sO_C(-m)^{p(m)})$ be the open set consisting
of locally free quotients $\sO_C(-m)^{p(m)}\to V'\to 0$ of rank
$r_1-r_2$ and degree $d_1-d_2$ such that $V'(m)$ is generated by
global sections, ${\rm H}^1(V'(m))=0$ and the quotient map induces
$\mathbb C^{p(m)}\cong {\rm H}^0(V'(m))$. Let $\sF\to \sH_2^0$ be
the frame bundle of $\pi_*(\sV'(m))$, where
$\pi:C\times\sH_2^0\to\sH_2^0$, then the pullback of the exact
sequence gives a morphism from $\sF$ to the projective bundle over
$Q\times\sU_2$ that parametrizes nontrivial extensions. The fiber
of this morphism has dimension at most $1$ since $V_1$ is a stable
bundle. Note that the irreducible component of $Q$ containing
stable bundles has maximal dimension and sending any extension
$(0\to V'\to V_1\to V_2\to 0)$ to ${\rm det}(V_2)^2\otimes {\rm
det}(V')$ defines a surjective morphism to $J^d(C)$. Thus
$$\text{dim($\sH_2$)}\le
(r_1-r_2)^2(g-1)+1+r_2^2(g-1)+1+(r,d)+r_2(r_1-r_2)(g-1)-g$$ and
codimension of $h(\sH_2)\subset\sR$ is at least
$r_1r_2(g-1)+1-(r,d)$.

To estimate $h(\sH_1)$, by Lemma \ref{lem3.1} (2), there are two
cases: (1) $r_1=r_2$, $(V_1,V_2)\in h(\sH_1)$ satisfy $0\to V_1\to
V_2\to \,_{x_{n_1}}\mathbb
C^{n_1}\oplus\cdots\oplus\,_{x_{n_k}}\mathbb C^{n_k}\to 0$ for
some $x_{n_i}\in C$ and $\sum n_i=d_2-d_1$ (the locus of these
points has codimension at least $r_1^2(g-1)+1-(r,d)$), or (2)
$\text{min$\{r_1,r_2\}$}>\frac{r}{(r,d)}$, $(V_1, V_2)\in
h(\sH_1)$ where $V_1$, $V_2$ are nontrivial extensions $0\to
V_1'\to V_1\to V_{k-1}\to 0$,  $0\to V_{k-1}\to V_2\to V_2'\to 0$
such that $V_{k-1}$ is a bundle of rank $r_{k-1}=k\frac{r}{(r,d)}$
and degree $d_{k-1}=k\frac{d}{(r,d)}$, where $1\le
k<\text{min$\{\frac{r_1(r,d)}{r},\,\frac{r_2(r,d)}{r}\}$}$. The
locus of such points has codimension at least
$r_{k-1}(r-r_{k-1})(g-1)+1-\frac{2r_{k-1}}{r}(r,d).$ Note that the
function $$f(x)=x(r-x)(g-1)+1-\frac{2x}{r}(r,d)$$ is an increase
function for $x\le\frac{r}{2}-\frac{(r,d)}{r(g-1)}$,
$r_0:=\frac{r}{(r,d)}\le
r_{k-1}\le\frac{r}{2}-\frac{(r,d)}{r(g-1)}$, and $f(r_1)\le
r_1(r-r_1)(g-1)+1-(r,d)$, we get \eqref{3.2} when
$\text{min$\{r_1,r_2\}$}>\frac{r}{(r,d)}$. If
$\text{min$\{r_1,r_2\}$}\le\frac{r}{(r,d)}$, any morphism $V_1\to
V_2$ must be surjective when $r_1>r_2$ and injective when
$r_1=r_2$. Thus we get the inequality \eqref{3.3}. The same
estimates also hold clearly for $\sD$.
\end{proof}

\begin{cor}\label{cor3.3} If $(r,d)\le 2$ and
$\sD\cap\sR\neq\emptyset$, then $\sD\cap\sR$, $\sD$ are
Cohen-Macaulay closed subschemes of codimension
$r_1r_2(g-1)+1-(r,d).$
\end{cor}

\begin{proof} By Lemma \ref{lem3.2}, when $(r,d)\le 2$, $\sD$ and
$\sD\cap\sR$ have codimension at least $r_1r_2(g-1)+1-(r,d)$. On
the other hand, it is standard to realize $\sD$ (resp.
$\sD\cap\sR$) as the degeneration locus of a morphism between two
vector bundles. Then the general theory implies that the
codimension of $\sD$ (resp. $\sD\cap\sR$) is at most
$r_1r_2(g-1)+1-(r,d)$ and $\sD$ (resp. $\sD\cap\sR$)
Cohen-Macaulay if the bound is reached.
\end{proof}

Write $\sP=P(r_1,d_1)$, $\sR=\sR(r_1,d_1)$. Recall that we have
$$\CD C\times\sP@>\pi>>\sP\\
@V1\times q VV  @V q VV\\
C\times\sR@>\pi>>\sR
\endCD$$ and the exact sequence
$$0\to (1\times q)^*\sV_1\otimes\pi^*\sO_{\sP}(1)\to\sE\to (1\times q)^*\sV_2
\to 0$$ which induces the morphism $$\Phi:\sP\to M.$$ Let
$Ad(\sE)$ denote the sheaf of trace free endomorphisms of $\sE$
and $\Delta(\sE)\subset Ad(\sE)$ the subsheaf of endomorphisms
that preserve the above exact sequence. Then
\ga{3.4}{0\to\Delta(\sE)\to Ad(\sE)\to(1\times
q)^*(\sV_1^{\vee}\otimes\sV_2)\otimes\pi^*\sO_{\sP}(-1)\to 0.} By
Lemma \ref{lem3.2},  $\pi_*(\sV_1^{\vee}\otimes\sV_2)=0$, thus the
sequence \eqref{3.4} induces \ga{3.5} {0\to R^1\pi_*\Delta(\sE)\to
R^1\pi_*Ad(\sE)\to
\\q^*R^1\pi_*(\sV_1^{\vee}\otimes\sV_2)\otimes\sO_{\sP}(-1)\to
0.\notag}

\begin{lem}\label{lem3.4} The infinitesimal deformation map $T_{\sP}\to
R^1\pi_*\Delta(\sE)$ induces an isomorphism. Under this
identification, the sequence \eqref{3.5} induces \ga{3.6}{0\to
T_{\sP}\xrightarrow{d\Phi}\Phi^*T_M\to
q^*R^1\pi_*(\sV_1^{\vee}\otimes\sV_2)\otimes\sO_{\sP}(-1)\to 0.}
\end{lem}

\begin{proof} Let $\sE nd^0={\rm ker}(\sE
nd(\sV_1)\oplus\sE nd(\sV_2)\xrightarrow{tr(\cdot)+tr(\cdot)}
\sO_{C\times\sR})$, then
$$0\to (1\times q)^*(\sV_1\otimes\sV_2^{\vee})\otimes\pi^*\sO_{\sP}(1)\to\Delta(\sE)\to (1\times
q)^*\sE nd^0 \to 0.$$ Now the proof is a straightforward
generalization of Lemma 6.6 in \cite{NR} since we have here
$T_{\sR}=R^1\pi_*\sE nd^0$, thus we omit it.
\end{proof}

For any $\xi=(V_1,V_2)\in\sR$, in order to study differential
$d\Phi_{\xi}$ of the morphism
$$\Phi_{\xi}:=\Phi|_{\sP_{\xi}}:\sP_{\xi}=q^{-1}(\xi)\to M,$$
let $[e]\in\sP_{\xi}$ be represented by a nontrivial extension
$$0\to V_1\xrightarrow{i}V\xrightarrow{j}V_2\to 0$$ and
$$K_{[e]}:=\{(f,g)\in\text{Hom}(V_1,V)\times\text{Hom}(V,V_2)\,|\,g\cdot
i+j\cdot f=0\}.$$

\begin{lem}\label{lem3.5} The kernel of $(d\Phi_{\xi})_{[e]}:
T_{\sP_{\xi},[e]}\to T_{M,\Phi([e])}$ has dimension
$$dim(K_{[e]})-1.$$
In particular, when $rk(V)=2$, $d\Phi_{\xi}$ is injective at every
point $[e]\in\sP_{\xi}$.
\end{lem}

\begin{proof} The $k[\e]$-value points of $\sP_{\xi}$ over
$[e]\in\sP_{\xi}$, which lie in kernel of $(d\Phi_{\xi})_{[e]}$,
are precisely represented by the extensions ($\e^2=0$)
$$0\to V_1\otimes_kk[\e]\xrightarrow{i_{\e}}V\otimes_kk[\e]\xrightarrow{j_{\e}}V_2\otimes_kk[\e]\to
0$$ with $i_{\e}=i\otimes 1+\e f\otimes 1$ and $j_{\e}=j\otimes
1+\e g\otimes 1$ where $(f,g)\in K_{[e]}$. Thus the kernel of
$(d\Phi_{\xi})_{[e]}$ has dimension $dim(K_{[e]})-1$. When
$rk(V)=2$, using Lemma \ref{lem3.1} (1), we can show Hom$(V_1,V)$
has dimension $1$, which implies the injectivity of
$(d\Phi_{\xi})_{[e]}$ (which also implies that $\Phi_{\xi}$ is an
embedding in the case of rank two).
\end{proof}

\begin{prop}\label{prop3.6} For any $\xi\in\sR\setminus \sR\cap\sD$,
the morphism $\Phi_{\xi}:\sP_{\xi}\to M$ is an embedding. For any
two different points $\xi_1,\,\xi_2\in\sR$, the intersection of
$\Phi_{\xi_1}(\sP_{\xi_1})$ and $\Phi_{\xi_2}(\sP_{\xi_2})$ has
dimension zero, i.e., a finite set.
\end{prop}

\begin{proof} $\xi=(V_1,V_2)\notin\sR\cap\sD$ means
$\text{Hom}(V_1,V_2)=0$, which implies that both $\Phi_{\xi}$ and
$d\Phi_{\xi}$ are injective, thus $\Phi_{\xi}$ is an embedding.

Let $\xi_1=(V_1,V_2)\in\sR$, $\xi_2=(W_1,W_2)\in \sR$ be any two
different points. Fix isomorphisms $\mathbb
P\cong\sP_{\xi_1}\cong\sP_{\xi_2}$ and pull back the universal
extensions to $C\times\mathbb P$
$$0\to p_1^*V_1\otimes\pi^*\sO_{\mathbb P}(1)\to\sE_1\to p_1^*V_2\to
0,$$ $$0\to p_1^*W_1\otimes\pi^*\sO_{\mathbb P}(1)\to\sE_2\to
p_1^*W_2\to 0$$ where $p_1:C\times\mathbb P\to C$,
$\pi:C\times\mathbb P\to\mathbb P$ are the projections. If the
intersection
$\Phi_{\xi_1}(\sP_{\xi_1})\cap\Phi_{\xi_2}(\sP_{\xi_2})$ has
positive dimension, then there is a nonsingular projective curve
$Y\to\mathbb P$ such that on $C\times Y$ the pullback of above
exact sequences
$$0\to p_1^*V_1\otimes\pi^*\sO_Y(1)\to\sE_1\to p_1^*V_2\to
0,$$ $$ 0\to p_1^*W_1\otimes\pi^*\sO_Y(1)\to\sE_2\to p_1^*W_2\to
0$$ define the same morphism $Y\to M$. Thus there is a line bundle
$\sN$ on $Y$ such that $\sE_1\cong \sE_2\otimes\pi^*\sN$. If
$\text{deg}(\sN)\le 0$, then
$$\text{Hom}(p_1^*V_1\otimes\pi^*\sO_Y(1),p_1^*W_2\otimes\pi^*\sN)=0$$
and $\sE_1\cong \sE_2\otimes\pi^*\sN$ induces an injection
$$p_1^*V_1\otimes\pi^*\sO_Y(1)\to p_1^*W_1\otimes\pi^*\sO_Y(1)\otimes\pi^*\sN,$$ which implies an
injection $V_1\to W_1\otimes\text{H}^0(\sN)$. Thus $\sN=\sO_Y$,
$V_1\cong W_1$ and $V_2\cong W_2$, which contradicts with
$\xi_1\neq\xi_2$. If $\text{deg}(\sN)\ge 0$, using
$\sE_2\cong\sE_1\otimes\pi^*\sN^{-1}$, we get contradiction by the
same arguments. Hence the intersection
$\Phi_{\xi_1}(\sP_{\xi_1})\cap\Phi_{\xi_2}(\sP_{\xi_2})$ has
dimension zero.
\end{proof}

It would be interesting to have a formula of the intersection
number of $\Phi_{\xi_1}(\sP_{\xi_1})$ and
$\Phi_{\xi_2}(\sP_{\xi_2})$.  We end this section with a question.

\begin{question} Is it true that any two lines on $M$ has at most
one intersection point ? It is interesting to describe the
configurations of lines on $M$ and on subvarieties (such as the
Brill-Noether locus) of $M$.\end{question}

\section{Remarks on minimal rational curves on $M$}

Let $M$ be the moduli space of stable bundles of rank $r$ and
degree $d$ with fixed determinant $\sL$ on a nonsingular
projective $C$ of genus $g\ge 3$. We assume $(r,d)=1$ in this
section. Then $M$ is a smooth projective Fano variety and there is
an universal bundle $\sE$ on $C\times M$. The universal bundle
$\sE$ is unique up to tensoring the pullback of a line bundle on
$M$. Since ${\rm Pic}(M)\cong\mathbb Z$, according to \cite[Remark
2.9]{Ra}, there is a unique universal bundle $\sE$ on $C\times M$
such that ${\rm det}(\sE|_{\{x\}\times M})=\Theta_M^{\alpha}$ for
any $x\in C$, where $\Theta_M$ is the ample generator of ${\rm
Pic}(M)$ and $\alpha$ is the smallest positive integer such that
$\alpha d\equiv\,1\,{\rm mod}\,(r)$. We will denote this canonical
universal bundle by $\sE$ in this section.

For any rational curve $\phi:\mathbb P^1\to M$ through a general
point of $M$, denote by $f:X=C\times\mathbb P^1\to C$ and
$\pi:X=C\times\mathbb P^1\to\mathbb P^1$ the projections, the
proof of \cite[Theorem 1]{Sun} implies in fact the following

\begin{thm}\label{thm4.1} If $\phi:\mathbb P^1\to M$ is a minimal
rational curve through a general point, then ${\rm
deg}(\phi^*\Theta_M)=r$ and $E:=(1\times\phi)^*\sE$ is a stable
bundle on $C\times\mathbb P^1$ with respect to any polarization.
Moreover, there is a point $x_{\phi}\in C$ such that
$E|_{\{x\}\times\mathbb P^1}=\sO_{\mathbb P^1}(\alpha)^{\oplus r}$
for $x\neq x_{\phi}$ and
$$E|_{\{x_{\phi}\}\times\mathbb P^1}=\sO_{\mathbb P^1}(\alpha+1)\oplus\sO_{\mathbb P^1}(\alpha)^{\oplus(r-2)}
\oplus\sO_{\mathbb P^1}(\alpha-1).$$ There is a stable vector
bundle $V$ on $C$ such that \ga{4.1}{0\to
f^*V\otimes\pi^*\sO_{\mathbb P^1}(\alpha)\to
E\to\sO_{\{x_{\phi}\}\times \mathbb P^1}(\alpha-1)\to 0} is an
exact sequence.
\end{thm}

For any general point $[W]\in M$, let $\Omega_W$ be the relative
cotangent bundle of $\mathbb P(W^{\vee})\to C$. Then \cite[Theorem
1]{Sun} also implies that the variety of all minimal rational
curves passing through $[W]\in M$ is naturally isomorphic to the
(double)projective bundle
$$\mathbb P(\Omega_W)\xrightarrow{p} C.$$ Thus, for
any $x_0\in C$, the set of minimal rational curves $\phi:\mathbb
P^1\to M$ with $x_{\phi}\neq x_0$ is the dense open set
$p^{-1}(C\setminus\{x_0\})$ of the variety of minimal rational
curves passing through $[W]\in M$. Let $\pi:C\times M\to M$ be the
projection and $d>2r(g-1)$. Then the direct image $\pi_*\sE$ is a
vector bundle on $M$ (called a Picard bundle). By using [9,
Theorem 1], we give a simple proof of some known results (See
\cite{BPN}, \cite{BBN} and \cite{LN}).

\begin{thm}\label{thm4.2} The bundles
$\sE_x:=\sE|_{\{x\}\times M}$ $(\forall x\in C)$, $\sE$ and the
Picard bundle $\pi_*\sE$ are stable with respect to any
polarization on $C\times M$ and $M$. Moreover, for any $x\neq y$,
we have $\sE_x\ncong\sE_y$.
\end{thm}

\begin{proof} By Proposition 3.7 in Chapter II of \cite{Ko}
(cf. also \cite{HM}, (4.3), proof of Proposition 12), for any closed
subset $S\subset M$ of codimension at least two, there is a minimal
rational curve $\phi:\mathbb P^1\to M$ such that $\phi(\mathbb
P^1)\cap S=\emptyset$. If $\sF\subset\sE_x$ is a subsheaf with
$\mu(\sF)\ge\mu(\sE_x)$, we may assume that the singular locus
$S\subset M$ of $\sF$ has codimension at least two. Then there is a
minimal rational curve $\phi:\mathbb P^1\to M$ with $x_{\phi}\neq x$
such that $\phi(\mathbb P^1)\cap S=\emptyset$ and $\phi(\mathbb
P^1)$ is not contained in the singular locus of $\sE_x/\sF$. By
Theorem \ref{thm4.1}, $E_x=\phi^*\sE_x=\sO_{\mathbb
P^1}(\alpha)^{\oplus r}$, thus $r\cdot a_{\sF}={\rm
deg}(\phi^*c_1(\sF))\le {\rm rk}(\sF)\cdot\alpha$ where
$a_{\sF}\in\mathbb Z$ such that $c_1(\sF)=a_{\sF}c_1(\Theta_M)$ and
$c_1(\sE_x)=\alpha c_1(\Theta_M)$, which implies
$\mu(\sF)=\mu(\sE_x)$, a contradiction since $\alpha
d\equiv\,1\,{\rm mod}\,(r)$. Thus $\sE_x$ ($\forall x\in C$) are
stable bundles.

To show stability of $\sE$ with respect to any polarization
$H=a\,f^{-1}(x)+b\,\Theta_M$, for any subsheaf $\sF\subset\sE$, let
$c_1(\sF)=d_1\, f^{-1}(x)+\beta\,\Theta_M$ and $c_1(\sE)=d\,
f^{-1}(x)+\alpha\, \Theta_M$, we have
$$\aligned&c_1(\sF)\cdot
H^n=(d_1b+\beta\,a)b^{n-1}f^{-1}(x)\cdot\Theta_M^n\\& c_1(\sE)\cdot
H^n=(d\,b+\alpha\,a)b^{n-1}f^{-1}(x)\cdot\Theta_M^n\endaligned$$
where $n={\rm dim}(M)$. Thus it is enough to show
\ga{4.2}{\frac{d_1b+\beta\,a}{{\rm
rk}(\sF)}<\frac{d\,b+\alpha\,a}{{\rm rk}(\sE)}.} We can assume that
singular loci $S\subset C\times M$ of $\sF$ has codimension at least
two. If $f(S)\subsetneqq C$, then stability of $\sE_x$ ($x\notin
f(S)$) and $\sE|_{C\times\{y\}}$ ($y\notin \pi(S)$) implies the
inequality \eqref{4.2}. If $f(S)=C$, for generic $x\in C$, the locus
$S_x=S\cap f^{-1}(x)\subset \{x\}\times M$ has codimension at least
two. Thus there is a minimal rational curve $\phi:\mathbb P^1\to M$
such that $\phi(\mathbb P^1)\cap \pi(S_x)=\emptyset$ and
$x_{\phi}\neq x$. Then $\{x\}\times\phi(\mathbb P^1)\subset C\times
M$ is disjoint with $S$ and $\sE|_{\{x\}\times\phi(\mathbb P^1)}$ is
semi-stable, which implies $\beta/{\rm rk}(\sF)\le \alpha/{\rm
rk}(\sE)$. The stability of $\sE|_{C\times\{y\}}$ ($y\notin \pi(S)$)
implies $d_1/{\rm rk}(\sF)<d/{\rm rk}(\sE)$. All together, we have
the inequality \eqref{4.2}.

To show stability of $\pi_*\sE$, for any subsheaf
$\sF\subset\pi_*\sE$, it is enough to find a $\phi:\mathbb P^1\to M$
disjoint with the singular locus $S$ of $\sF$ such that the
restrictions $F=\phi^*\sF$ and $\pi_*E=\pi^*(\pi_*\sE)$ satisfy
$\mu(F)<\mu(\pi_*E)$, where $E=(1\times\phi)^*\sE$. Let
$\sF(W)\subset{\rm H}^0(W)=\pi_*(\sE)|_{[W]}$ be the fibre of $\sF$
at a general point $[W]\in M$. Let $Z\subset C$ be the set of common
zero points of sections of $\sF(W)$ and $x\in C\setminus Z$ a
general point. Let $\sF(W)_x=\{s_x\in W_x\,|\,s\in\sF(W)\}$ and
$\zeta\in \mathbb P(W^{\vee}_x)$ a general point such that
$\sF(W)_x\nsubseteq\zeta^{\bot}\subset W_x$. Define a vector bundle
$W^{\zeta}$, which is the Hecke modification of $W$ along
$\zeta^{\bot}\subset W_x$, by
$$0\to W^{\zeta}\xrightarrow\iota W\to (W_x/\zeta^{\bot})\otimes\sO_x\to 0$$
where $\zeta^{\bot}$ denotes the hyperplane in $W_x$ annihilated
by $\zeta$. The $1$-dimensional subspace ${\rm
ker}(\iota_x)\subset W_x^{\zeta}$ defines a point $[{\rm
ker}(\iota_x)]\in \mathbb P(W_x^{\zeta})$. Then a general line
$\ell\subset\mathbb P(W_x^{\zeta})$ passing through $[{\rm
ker}(\iota_x)]\in \mathbb P(W_x^{\zeta})$ defines a minimal
rational curve $\phi:\mathbb P^1\to M$ passing through $[W]\in M$
disjoint with $S$ such that $x_{\phi}=x$. By \eqref{4.1}, we have
\ga{4.4}{0\to {\rm H}^0(V)\otimes\sO_{\mathbb P^1}(\alpha)\to
\pi_*E\to\sO_{\{x_{\phi}\}\times \mathbb P^1}(\alpha-1)\to 0.}
Since $\sF(W)_x\nsubseteq\zeta^{\bot}\subset W_x$ and
$\phi(\mathbb P^1)$ passes through $[W]\in M$, the image of
$F\subset \pi_*E$ under the surjection
$\pi_*E\to\sO_{\{x_{\phi}\}\times \mathbb P^1}(\alpha-1)$ is
non-trivial. Thus
$$\mu(F)\le \alpha -\frac{1}{{\rm rk}(F)}<\alpha-\frac{1}{{\rm
rk}(\pi_*E)}=\mu(\pi_*E).$$

To show $\sE_x\ncong \sE_y$ when $x\neq y$, we choose a minimal
rational curve $\phi:\mathbb P^1\to M$ with $x_{\phi}=x$. Then, by
Theorem \ref{4.1}, we have
$$\phi^*\sE_y=\sO_{\mathbb P^1}(\alpha)^{\oplus r}\neq\sO_{\mathbb P^1}(\alpha+1)\oplus\sO_{\mathbb P^1}(\alpha)^{\oplus(r-2)}
\oplus\sO_{\mathbb P^1}(\alpha-1)=\phi^*\sE_x.$$ Thus
$\sE_x\neq\sE_y$, we finish the proof of theorem.
\end{proof}

\begin{rmk}\label{rmk4.3} As far as we know, the semi-stability of $\sE_x$
appears first of all as Proposition 1.4 in \cite{BPN}, its
stability is Proposition 2.1 of \cite{LN}. The stability of $\sE$
is Theorem 1.5 of \cite{BPN}. The stability of $\pi_*\sE$ and the
fact that $\sE_x\ncong\sE_y$ ($x\neq y$) are the main theorems of
\cite{BBN} and \cite{LN}.\end{rmk}

\bibliographystyle{plain}

\renewcommand\refname{References}

\end{document}